\documentclass[a4paper]{article}

\usepackage[english]{babel}
\usepackage[utf8x]{inputenc}
\usepackage[T1]{fontenc}

\usepackage[a4paper,top=3cm,bottom=2cm,left=3cm,right=3cm,marginparwidth=1.75cm]{geometry}

\usepackage{booktabs}
\usepackage{amsmath}
\usepackage{amssymb}
\usepackage{graphicx}
\usepackage[colorinlistoftodos]{todonotes}
\usepackage[colorlinks=true, allcolors=blue]{hyperref}
\usepackage{url}
\usepackage{caption}
\DeclareCaptionType{Figure}

\newcommand{\asvector}[1]{{#1}}
\newcommand{\tangentspace}[1]{\mathrm{T}_{\!#1}}
\newcommand{\gradient}[3]{\nabla_{#1}^{#2}\!#3}
\newcommand{\maxsubscript}[1]{#1_{\mathrm{max}}}
\newcommand{\projection}[0]{\rho}
\newcommand{\reals}[0]{\mathbb{R}}
\newcommand{\area}[1]{\mathrm{A} (#1)}
\newcommand{\frechet}[1]{\mu (#1)}
\newcommand{\collection}[1]{\mathbf #1}
\newcommand{\sphere}[1]{\mathbb{S}^{#1}}
\newcommand{\hyperboloid}[1]{\mathbb{H}^{#1}}
\newcommand{\pb}[1]{\mathbb{B}^{#1}}  
\newcommand{\minkowski}[1]{\mathbb{R}^{#1:1}}
\newcommand{\form}[2]{\left\langle {#1}, {#2} \right\rangle}
\newcommand{\minkowskiformpower}[4]{\left\langle {#2},{#3} \right\rangle_{#1:1}^{#4}}
\newcommand{\minkowskiform}[3]{\minkowskiformpower {#1} {#2} {#3} {}}
\newcommand{\new}[1]{{#1}^{\mathrm{new}}}
\newcommand{\dist}[0]{\mathrm{d}}
\newcommand*{\differential}{\mathop{}\!\mathrm{d}}

\DeclareMathOperator{\Exp}{Exp}

\DeclareMathOperator*{\argmin}{argmin}
\DeclareMathOperator{\arccosh}{arccosh}

\title{Gradient descent in hyperbolic space}

\author{
  Benjamin Wilson\\
  Lateral GmbH\\
  \texttt{benjamin@lateral.io}
  \and
  Matthias Leimeister\\
  Lateral GmbH\\
  \texttt{matthias@lateral.io}
}

\begin{document}
\maketitle

\begin{abstract}
Gradient descent generalises naturally to Riemannian manifolds \cite{Bonnabel}, and to hyperbolic $n$-space, in particular.  Namely, having calculated the gradient at the point on the manifold representing the model parameters, the updated point is obtained by travelling along the geodesic passing in the direction of the gradient.
Some recent works employing optimisation in hyperbolic space have not attempted this procedure, however, employing instead various approximations to avoid a calculation that was considered to be too complicated.
In this tutorial, we demonstrate that in the hyperboloid model of hyperbolic space, the necessary calculations to perform gradient descent are in fact straight-forward.  The advantages of the approach are then both illustrated and quantified for the optimisation problem of computing the Fréchet mean (i.e.\ barycentre) of points in hyperbolic space.
\end{abstract}

\begin{section}{Introduction}
Hyperbolic $n$-space is an $n$-dimensional Riemannian manifold with a constant sectional curvature of $-1$.
Like Euclidean space, it is unbounded in its extent and homogeneous (``isotropic'').
Unlike Euclidean space, the circumference of a circle grows exponentially in the radius.
This property makes hyperbolic space particularly well-suited for the embedding of trees, since the rate of growth of the circumference is able to match the rate of growth of the number of nodes as the depth of the tree increases.  Indeed, the superiority of hyperbolic space for the embedding of tree-like graphs has been dramatically demonstrated in the recent works \cite{NickelKiela, DeSa}\footnote{The recent graph embedding of \cite{Chamberlain} is not (despite its title) in hyperbolic space, since the optimisation takes place in the tangent space of a single point, which is Euclidean.}.
The related problem of multi-dimensional scaling in hyperbolic space has also been approached in \cite{LindmanCaelli,WilsonHancock, DeSa, CvetkovskiCrovella2016}.

Gradient descent generalises naturally to Riemannian manifolds \cite{Bonnabel}, and to hyperbolic $n$-space, in particular.
Having calculated the gradient at the point on the manifold representing the model parameters, gradient descent updates the point by travelling along the geodesic that passes through the point in the opposite direction to the gradient.
In other words, the ``exponential map'' is applied to an appropriate multiple of the gradient vector.
Moreover, if the hyperboloid model of hyperbolic space is used, then necessary calculations are straight-forward.

However, despite the straight-forwardness of gradient descent in the hyperboloid model, the above-mentioned optimisations have employed more complicated formulations that in some cases are merely approximations of gradient descent.
Instead of using gradient descent on hyperbolic space, \cite{WilsonHancock} perform the optimisation entirely in the tangent space of the current point, using the exponential map only to map the convergent back onto the manifold.  As the tangent space is not isometric with the manifold, this approach ignores important aspects of the geometry.
Earlier work \cite{LindmanCaelli} approached the same problem, but used a ``direct search'' method instead of gradient descent, in order to avoid perceived complexities.
The works \cite{NickelKiela,DeSa}, on the other hand, worked in the Poincaré ball model of hyperbolic space, where additive updates are possible.  These additive updates approximate the exponential map and thereby also gradient descent.  While this approach has theoretic underpinning in the notion of a ``retraction'' \cite{Bonnabel}, it introduces significant (and as we'll see, in the case of hyperbolic space, unnecessary) inaccuracy into the optimisation.  It suffers further from the drawback that the gradient of the distance function on the Poincaré ball is complicated to compute.
Finally, \cite{CvetkovskiCrovella2016} perform a steepest descent line search on the Poincaré disc using Möbius transformations, which is applicable only in the 2-dimensional case.

In this work, we describe the hyperboloid model of hyperbolic space  and demonstrate its advantages for gradient optimisation. Indeed, it is shown that optimisation using the hyperboloid model is precisely as straight-forward as optimisation on the sphere.
Specifically, it is shown that the gradient of any differentiable objective function written in terms of the distance on the manifold is trivial to compute and moreover that the gradient update via the exponential map can be calculated using a simple combination of the hyperbolic functions $\cosh$ and $\sinh$\footnote{Recall that $\cosh r = \frac{1}{2}(e^r + e^{-r})$ and $\sinh r = \frac{1}{2}(e^r - e^{-r})$. }.
The hyperboloid model is well-known in many fields, e.g.\ image processing \cite{Bergmann2017}. The facilities it offers for gradient optimisation seem to be under-appreciated by some other machine learning practitioners, however.
This is evidenced by the use of the (computationally much more challenging, yet equivalent) Poincaré ball model in \cite{NickelKiela,DeSa}.
That it is not fully-understood can be further seen by the optimisations of \cite{LindmanCaelli,WilsonHancock}, which \textit{do} employ the hyperboloid model, but avoid the calculations needed for gradient descent.
The objective of this work is to demonstrate that gradient descent using the hyperboloid model can indeed be very simple.

Section \ref{section:optimisationonthesphere} reviews well-known facts about optimisation on the sphere.
This greatly facilitates the introduction of the hyperboloid model of hyperbolic space in section \ref{section:optimisationonthehyperboloid}.
Indeed, optimisation on the hyperboloid is closely analogous to optimisation on the sphere.
In section \ref{section:thepoincaremodel}, the Poincaré model of hyperbolic space is introduced and its relationship to the hyperboloid model is discussed.
This includes a discussion of the ``retraction'' update, used by \cite{NickelKiela,DeSa} to approximate the exponential map and a measurement of its imprecision.
Finally, in section \ref{section:numericalexperiments}, the advantages of the hyperboloid model are then both illustrated and quantified in the case of the optimisation problem of computing the Fréchet mean (the analogue of the barycentre) of points in hyperbolic space.
Specifically, in this case, it is shown that the exponential updates on the hyperboloid arrive in the neighbourhood of the solution in 46\% less gradient updates (on average) than using the ``retraction'' updates on the Poincaré ball.
\end{section}

\begin{section}{Recent works}
	The exponential map on the hyperboloid (transported to the Poincaré ball) is used for full Riemannian gradient descent in the recent work \cite{Ganea2018}.
	Furthermore, two related works have appeared since the publication of the first version of this paper.
	In \cite{NickelKiela2018}, full Riemannian gradient descent in hyperbolic space is employed for embedding word graphs.
	The derivation of the algorithm there is thus largely similar to our presentation.
	In \cite{Enokida2018}, on the other hand, a formula for the exponential map on the Poincaré disc model of 2-dimensional hyperbolic space is derived directly.  As we shall see below, this formula finds a much simpler expression in the hyperboloid model (where it is valid in any dimension).
\end{section}

\begin{section}{A review of optimisation on the sphere}\label{section:optimisationonthesphere}
The relationship of the hyperboloid to its ambient is closely analogous to the relationship of the sphere and its Euclidean ambient, and it is instructive to consider this case first.
Let $n \geqslant 1$ and write
$$ \sphere n = \{\, \asvector x \in \reals^{n+1} \mid \| \asvector x \| = 1 \,\} $$
for the $n$-dimensional sphere $\sphere n$ with unit radius embedded in $(n+1)$-dimensional Euclidean space $\reals^{n+1}$.

Recall that the distance between two points on $\sphere n$ is just their planar angle in the ambient, and so can be computed via
\begin{equation}\label{sphericaldistance}
	\dist_{\sphere n} (\asvector u, \asvector v) = \arccos(\form {\asvector u} {\asvector v}), \quad \asvector u , \asvector v \in \sphere n,
\end{equation}
and furthermore that, at each point $\asvector p \in \sphere n$ there is an $n$-dimensional tangent space
$$\tangentspace {\asvector p} \sphere n = \{\, {\asvector x} \in \reals^{n+1} \mid \form {\asvector p} {\asvector x} = 0 \,\}$$
consisting of all those vectors from the ambient space $\reals^{n+1}$ that are perpendicular to ${\asvector p}$.
When optimising ${\asvector p}$ with respect to some error function $E$ defined on the sphere, the gradient vector $\gradient {{\asvector p}} {\sphere n} E \in \tangentspace {\asvector p} \sphere n$ is the vector in the tangent space that points in the direction of greatest (instantaneous) increase in $E$, and its length gives the partial derivative of $E$ in that direction.
In order to perform a gradient update, the point $\asvector p$ is moved distance $\| \gradient p {\sphere n} E \|$ along the great circle passing through $\asvector p$ in the direction of the negative of the gradient vector.
This new point is called the {\it exponential} $\Exp_{\asvector p} (-\gradient p {\sphere n} E)$ of the negative gradient at p.
Indeed, in the case of the sphere, the exponential map $\Exp_{\asvector p}$ is defined on the entire tangent space
$$ \Exp_{\asvector p} : \tangentspace {\asvector p} \sphere n \to \sphere n, \quad \asvector v \mapsto \Exp_{\asvector p}(\asvector v)$$
and is given by the convenient formula
$$ \Exp_{\asvector p}(\asvector v) = \cos(\| \asvector v \|) \asvector p + \sin(\| \asvector v \|) \frac{\asvector v}{\| \asvector v \|}.$$
Notice that this formula is in terms of the vector arithmetic and inner product of the ambient linear space $\reals^{n+1}$.  This is representative of the very helpful relationship between the geometry of the sphere and the arithmetic of its ambient space.
If the error function $E$ is computable near the sphere as well as on it (as is the case with the distance formula \eqref{sphericaldistance}, for example), then the ambient can help us again.
For by considering $\asvector p$ as a point of the ambient space, the gradient $\gradient p {\reals^{n+1}} E$ at $\asvector p$ can be easily calculated as the vector of partial derivatives at $\asvector p$.
Using the dot product, this vector can then be projected onto $\tangentspace p \sphere n $, the tangent space of the sphere at $\asvector p$, to obtain the gradient on the sphere:
$$ \gradient {\asvector p} {\sphere n} E = \gradient {\asvector p} {\reals^{n+1}} E - \form {\asvector p} {\gradient {\asvector p} {\reals^{n+1}} E} \cdot \asvector p .$$
Using the ambient thus allows the computation of the gradient without recourse to a choice of coordinates on the sphere.
\end{section}

\begin{section}{The hyperboloid model of hyperbolic space}\label{section:optimisationonthehyperboloid}
Hyperbolic space can not be embedded without distortion in Euclidean space \cite{Efimov}, but there are various {\it models} of hyperbolic space that allow calculations to be carried out.
As we'll see, the hyperboloid model is particularly well-suited to deriving formulae e.g.\ for optimisation.
A broader discussion of the properties of the hyperboloid model is offered by \cite{Reynolds}.

\begin{subsection}{The hyperboloid model}
The unit sphere $\sphere n$ is a space of constant positive curvature.  Hyperbolic space has constant {\it negative} curvature, and can be constructed, analogously to the sphere in its Euclidean ambient, as a pseudo-sphere (or hyperboloid) in a linear ambient space called {\it Minkowski space}.
This construction is called the ``hyperboloid model'' of hyperbolic space.
It is well-known in both mathematics and physics and has been used for optimisation in \cite{WilsonHancock,LindmanCaelli}.
However, even these papers did not fully exploit the helpful relationship (analogous to the spherical case, above) between the hyperboloid and its ambient vector space.

Write $\minkowski n$ for a copy of $\reals^{n+1}$ equipped with the bilinear form $\minkowskiform n \cdot \cdot$, given by:
$$ \minkowskiform n {\asvector u} {\asvector v} = \sum_{i=1}^n u_i v_i - u_{n+1} v_{n+1}.$$
This is the $(n+1)$-dimensional {\it Minkowski space}.
Notice that, in contrast to the familiar dot product on Euclidean space, the bilinear form of Minkowski space is not positive definite, i.e. there exist vectors $\asvector v$ such that $\minkowskiform n{\asvector v} {\asvector v} < 0$.
Indeed, the {\it $n$-dimensional hyperboloid} $\hyperboloid n$ is a collection of such points:
$$ \hyperboloid n = \{\, \asvector x \in \minkowski n \mid \minkowskiform n {\asvector x} {\asvector x} = -1,\ x_{n+1} > 0 \,\}. $$
The distance between two points on the hyperboloid can be easily computed via
\begin{equation}\label{hyperboloiddistance}
	\dist_{\hyperboloid n} (\asvector u, \asvector v) = \arccosh(- \minkowskiform n {\asvector u} {\asvector v}), \quad \asvector u, \asvector v \in \hyperboloid n,
\end{equation}
(compare \eqref{sphericaldistance}).
As in that case, the tangent space to a point $\asvector p \in \hyperboloid n$ is the set of all vectors perpendicular to $\asvector p$ with respect to the bilinear form of the ambient
$$\tangentspace p \hyperboloid n = \{\, x \in \minkowski n \mid \minkowskiform n {\asvector p} {\asvector x} = 0 \,\}.$$
Each tangent space inherits the bilinear form from the ambient, and this restriction is positive-definite (thus $\hyperboloid n$ is a Riemannian manifold, embedded in the {\it pseudo}-Riemannian ambient $\minkowski n$).
In particular, for any $\asvector p \in \hyperboloid n$, it makes sense to talk about the norm $\| \asvector v \|$ of a tangent vector $\asvector v \in \tangentspace {\asvector p} \hyperboloid n$.
\end{subsection}

\begin{subsection}{Optimisation on the hyperboloid $\hyperboloid n$}
Analogously to the spherical case, the gradient of an error function $E$ at a point $\asvector p$ is a vector in the tangent space, and gradient optimisation is achieved by applying the exponential map at $\asvector p$ to the gradient vector (or some multiple thereof).  For $\asvector p\in \hyperboloid n$, the exponential map is given by
\begin{equation}\label{hyperboloidexponential}
\Exp_{\asvector p} (v) = \cosh (\| \asvector v \|) \asvector p + \sinh (\| \asvector v \|) \frac{\asvector v}{\| \asvector v \|}.
\end{equation}
As in the case of the sphere, this formula is very simple.
Moreover, the formula allows the calculation of a gradient update precisely, without the need for any approximation (such as is given by the retraction on the Poincaré ball, discussed below).
If the error function is also defined on a region of $\minkowski n$ around $\hyperboloid n$ (as is the case e.g.\ for the distance function \eqref{hyperboloiddistance}), then the ambient helps with the calculation of the gradient as before, with two small modifications.
Firstly, the gradient $\gradient {\asvector p} {\minkowski n} E \in \tangentspace p \minkowski n$ in the ambient space is the vector of partial derivatives given by
$$\gradient {\asvector p} {\minkowski n} E = \textstyle \left( \frac{\partial E}{\partial x_1} |_p, \ldots, \frac{\partial E}{\partial x_{n}} |_p, -\frac{\partial E}{\partial x_{n+1}} |_p \right),$$
(note the minus sign).
For example, since
$$\frac{\partial}{\partial u_i} \minkowskiform n {\asvector u} {\asvector v} = (-1)^{\delta_{i, n+1}} v_i,$$
for any $1 \leqslant i \leqslant n+1$, the gradient (in the ambient) of the distance function \eqref{hyperboloiddistance}, with respect to one of its arguments $\asvector u$, has the simple form\footnote{Recall that $\frac{\mathrm d}{\mathrm{d}z} \arccosh(z) = (z^2 - 1)^{-1/2}$ }
\begin{equation}\label{hyperbolicdistancegradient}
	\gradient {\asvector u} {\minkowski n} {\dist_{\hyperboloid n} (\asvector u, \asvector v)} = -({\minkowskiformpower n {\asvector u} {\asvector v} 2} - 1)^{-1/2} \cdot \asvector v.
\end{equation}
The simplicity of this formula (contrast with the formula in \cite{NickelKiela}, see section \ref{subsection:retraction}) leads to simple expressions for the gradients of the error functions defined in terms of the distance (e.g.\ those used in \cite{NickelKiela,DeSa,WilsonHancock}, when re-expressed on the hyperboloid).
The second important difference from the spherical case is in the projection of a vector from the ambient onto the tangent space of the hyperboloid, which is given by:
\begin{equation}\label{hyperboloidtangentspaceprojection}
\gradient {\asvector p} {\hyperboloid n} E = \gradient {\asvector p} {\minkowski n} E + \minkowskiform n {\asvector p} {\gradient {\asvector p} {\minkowski n} E} \cdot \asvector p.
\end{equation}
(again, note the difference in sign).
Figure \ref{gradientdescentalgorithm} summarises the resulting algorithm for gradient optimisation.

\begin{center}
\fbox{
	\centering
    \parbox{\textwidth}{
Given an error function $E$ defined in terms of the distance, a learning rate $\alpha$ and an initial value $\Theta = \Theta^{(0)}$, repeat the following until convergence:
\begin{enumerate}
	\item Calculate the gradient $\gradient \Theta {\minkowski n} E$ using \eqref{hyperbolicdistancegradient}.
	\item Calculate $\gradient \Theta {\hyperboloid n} E$ from $\gradient \Theta {\minkowski n} E$, using \eqref{hyperboloidtangentspaceprojection}.
	\item Set $\new \Theta = \Exp_\Theta (-\alpha \cdot \gradient \Theta {\hyperboloid n} E)$, calculated using \eqref{hyperboloidexponential}.
\end{enumerate}
\captionof{Figure}{Algorithm for gradient descent on $\hyperboloid n$\label{gradientdescentalgorithm}.}
}}
\end{center}
\end{subsection}
\end{section}

\begin{subsection}{Choosing the learning rate}
As in the case of gradient descent in Euclidean space, a crucial aspect of the algorithm in Figure~\ref{gradientdescentalgorithm} is the choice of an appropriate learning rate (or ``step size'') $\alpha$.
In practical applications, instead of running until convergence, often early stopping after a pre-defined number of iterations is used and the learning rate is either constant or adjusted during the training by some schedule. For example, fixing the number of iterations and linearly decaying the learning rate to zero is one possibility. In order to guarantee convergence of gradient descent in hyperbolic space, several results on the appropriate choice of learning rates are available. The \textit{Armijo rule} (\cite{Absil}, Definition 4.2.2)  gives bounds on the learning rate in terms of the norm of the gradient of the objective function. The same work proposes an accelerated line search algorithm (\cite{Absil}, p. 63, Algorithm 1) that, when using the Armijo rule, converges to a local optimimum of the objective function. In another work, \cite{Afsari} investigate the convergence of Riemannian gradient descent with constant step sizes. The step sizes that guarantee convergence of gradient descent for the Fréchet mean depend on the Hessian of the objective function and the radius of a geodesically convex ball bounding the data points.
However, for simplicity, the numerical experiments of section \ref{section:numericalexperiments} use constant learning rates fixed independently of the data points.
\end{subsection}

\begin{section}{The Poincaré ball model and the retraction}\label{section:thepoincaremodel}
\begin{subsection}{The Poincaré ball model}
We saw in section \ref{section:optimisationonthehyperboloid} that the hyperboloid is a model of hyperbolic space that is particularly well-suited to the calculations required for gradient descent.
The Poincaré ball is another model of hyperbolic space.
Carrying out the same calculations using the Poincaré model (without using the hyperboloid!) is a complicated undertaking.
The strength of the Poincaré ball model is that it permits a form of visualisation of hyperbolic space and of any arrangement of points within it.
Fortunately, as we'll see, it is easy to move between the two models.

The Poincaré ball model represents $n$-dimensional hyperbolic space as the interior of the unit ball in $n$-dimensional Euclidean space:
$$ \pb n = \{\, \asvector x \in \reals^n \mid \| \asvector x \| < 1 \,\}.$$
In this model, the distance between two points can be calculated using the ambient Euclidean geometry via
\begin{equation}\label{pbdistance}
	\dist_{\pb n} (\asvector u, \asvector v) = \arccosh \left(1 + 2 \frac{\| \asvector u - \asvector v \|^2}{(1 - \| \asvector u \|^2)(1 - \| \asvector v \|^2)}\right)
\end{equation}
In the case where $n=2$, the Poincaré ball model permits a very convenient (if distorted) visualisation of hyperbolic space.
This distortion manifests itself in the dilation of distance between a pair of points by the closeness of either to the Euclidean boundary of the ball, and (consequently) in the seemingly distinguished nature of the centre point.
In fact, hyperbolic space is isotropic, meaning that there are no distinguished points.
\end{subsection}

\begin{subsection}{Relationship to the hyperboloid model}
It is easy to map points from the hyperboloid model to their correspondents in the Poincaré ball model, via:
\begin{equation}\label{hyperboloidtopb}
	\projection: \hyperboloid n \to \pb n, \quad \projection(\asvector x) = \frac{1}{x_{n+1} + 1}(x_1, \ldots, x_n), \quad \asvector x \in \hyperboloid n
\end{equation}
This is analogous to a map projection of the sphere $\sphere 2$.  The inverse of this map is given by:
\begin{equation}\label{pbtohyperboloid}
	\projection^{-1}: \pb n \to \hyperboloid n, \quad \projection^{-1}(\asvector y) = \frac{2}{1 - r^2}\left(y_1, \ldots, y_n, \frac{1 + r^2}{2}\right), \quad \asvector y \in \pb n,
\end{equation}
where $r = \| \asvector y \|$ is the Euclidean norm of $\asvector y$. 
The following formula for the expression of the differential $\differential \projection |_{\asvector x} : \tangentspace {\asvector x} \hyperboloid n \to \tangentspace {\projection({\asvector x})} \pb n$ allows us to map a hyperboloid gradient to the corresponding gradient on the Poincaré ball, and will be useful in our experiments.
For $\asvector x \in \hyperboloid n$ and $v \in \tangentspace {\asvector x} \hyperboloid n$, we have
\begin{equation}\label{differentialhyperboloidtopb}
	(\differential \projection |_{\asvector x} (v))_i = \frac{1}{x_{n+1} + 1} \left(v_i - \frac{x_i v_{n+1}}{x_{n+1} + 1}\right), \quad 1 \leqslant i \leqslant n.
\end{equation}
\end{subsection}

\begin{subsection}{Retraction as approximation of the exponential map}\label{subsection:retraction}
The exponential map, applied to a one-dimensional subspace of the tangent space, yields a circular line segment that meets the boundary of the Poincaré ball at right angles, or a diameter of the ball.
However, despite being easy to describe geometrically, the exponential map on the Poincaré ball is difficult to calculate.
To avoid this calculation, recent papers \cite{NickelKiela,DeSa} have used a first-order approximation to the exponential map, called a {\it retraction} (see \cite{Bonnabel}), which performs a gradient step by simply adding a scalar multiple of the gradient vector.
That is, for a parameter vector $\asvector \Psi \in \pb n$:
\begin{equation}\label{retraction}
	\new {\asvector \Psi} = \asvector \Psi - \alpha \cdot \gradient {\asvector \Psi} {\pb n} E
\end{equation}
where $E$ is the error function defined on $\pb n$ and $\alpha > 0$ is a learning rate.
The gradient $\gradient {\asvector \Psi} {\pb n} E$ can moreover be calculated from the Euclidean gradient of $E$ by a simple re-scaling.
While the simplicity of this approach is very attractive, it introduces errors in the gradient update that are very significant when the point ${\asvector \Psi}$ being updated is already some distance from the centre point of the ball (see Figure \ref{retractionerrorfigure}).
Indeed, for very large step sizes, it is even possible that $\new {\asvector \Psi} \notin \pb n$, i.e.\ that the retraction update erroneously skips over an infinite expanse of parameter space and then leaves the space entirely.
In this case, the practice is to re-scale $\new {\asvector \Psi}$ such that it is again in the parameter space \cite{NickelKiela}.

\begin{figure}
	\centering
	\includegraphics[width=9cm]{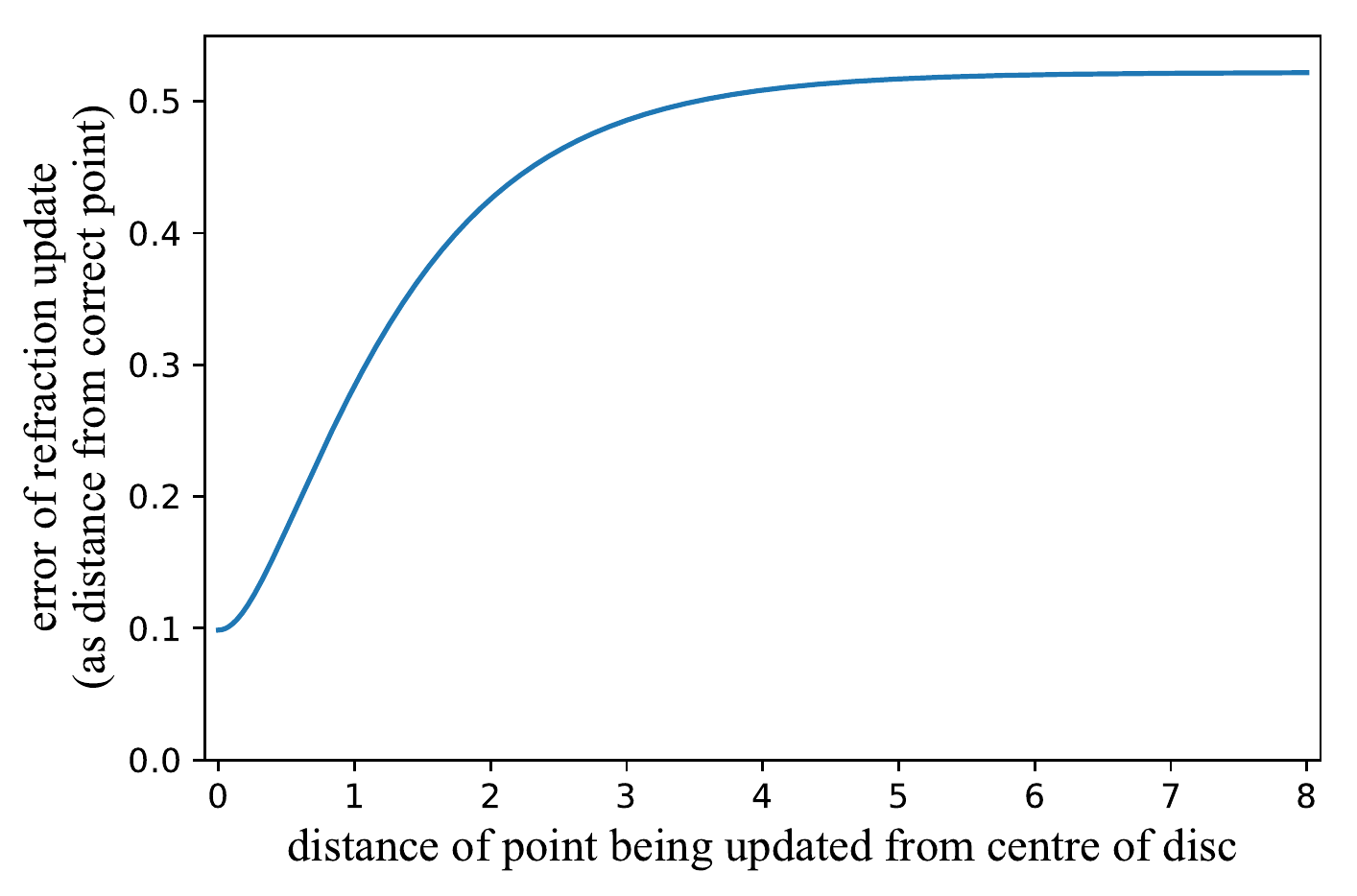}
	\captionof{Figure}{Worst-case inaccuracy of the retraction approximation for a unit length step.\label{retractionerrorfigure}}
\end{figure}

As we saw above, using the hyperboloid model of hyperbolic space permits the calculation of the precise update afforded by the exponential map in an elementary manner.
This avoids the unnecessary, high-error approximation of the retraction.
\end{subsection}

\begin{subsection}{Gradient of the distance function}
The second important computational disadvantage of the Poincaré ball model is the calculation of the gradient of the distance function.  The following expression is from the recent work \cite{NickelKiela}:
\begin{equation}
	\gradient {\asvector u} {\pb n} {\dist ({\asvector u}, {\asvector v})} = \frac{4}{b \sqrt{c^2 - 1}} \left( \frac{\| {\asvector v} \|^2 - 2 \form {\asvector u} {\asvector v} + 1}{a^2} {\asvector u} - \frac{{\asvector v}}{a} \right),
\end{equation}
where ${\asvector u}, v \in \pb n$, and
$$ a = 1 - \| {\asvector u} \|^2, \quad b = 1 - \| {\asvector v} \|^2, \quad c = 1 + \frac{2}{ab}\| {\asvector u} - {\asvector v}\|^2.$$
The complexity of this evpression is also present, of course, in the gradient on the Poincaré ball of any objective function that is a function of the distance.
The gradient of the distance function in the hyperboloid model is, in contrast, straight-forward to calculate.
\end{subsection}
\end{section}

\begin{section}{Numerical experiments}\label{section:numericalexperiments}
	In this section we demonstrate the advantage of using the exponential map \eqref{hyperboloidexponential} on the hyperboloid instead of the retraction \eqref{retraction} on the Poincaré disc by considering the problem of finding the barycentre (or centre of mass) of given points in hyperbolic space.
	As described below, this is formulated as an optimisation problem.
	For simplicity, we use constant learning rates chosen independently of the given points.
The experiment is depicted on the Poincaré disc in Figure \ref{frechetmeanfigure}.
The source code for the experiments is available online\footnote{\url{https://github.com/lateral/frechet-mean-hyperboloid}}.

\begin{figure}
	\centering
	\includegraphics[width=15cm]{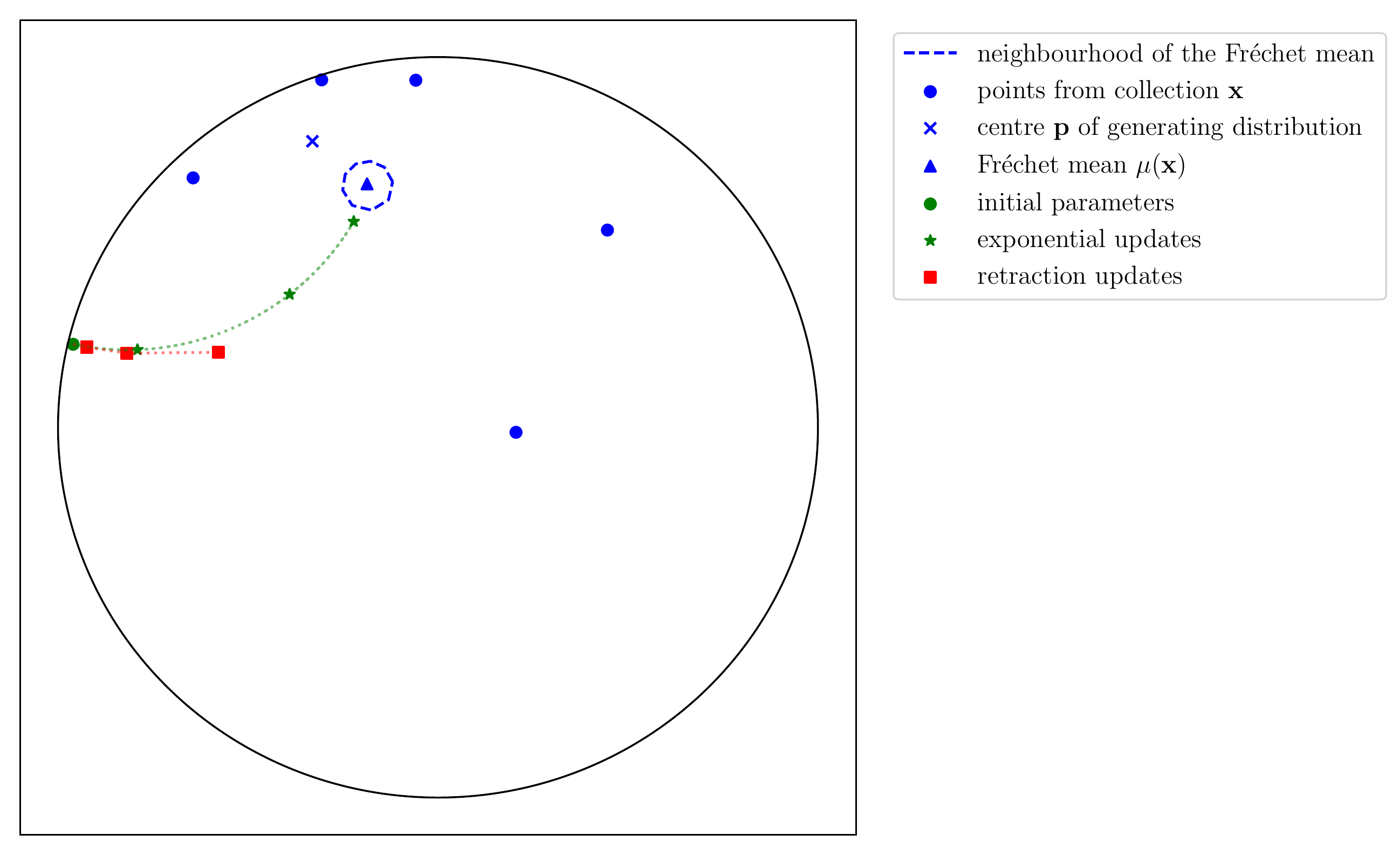}
	\captionof{Figure}{Fréchet mean optimisation experiment depicted on the Poincaré disc.
	Beginning from the same initial value and using the same learning rate (here $0.45$), exponential- and retraction- updates are computed.  The number of steps before arriving in neighbourhood about the Fréchet is counted in both cases.
Neighbourhood radius is here very large ($0.25$) for illustrative purposes.\label{frechetmeanfigure}}
\end{figure}

\begin{subsection}{The Fréchet mean}\label{section:numericalexperiments}
	Let $s$ be a positive integer, and denote by $\collection x = \{ {\asvector x}^{(1)}, \ldots, {\asvector x}^{(s)}\}$ a set of $s$ distinct points.
In Euclidean space, the barycentre (or centre of mass) $\frechet {\collection x}$ of the points $\collection x$ is simply their arithmetic average, i.e.\ $\frechet {\collection x} = \frac{1}{s} \sum_{i=1}^s {\asvector x^{(i)}}$.
This is not the case in any of the models of hyperbolic space, however.  In these cases (and for more general Riemannian manifolds) the problem of finding the barycentre must be posed instead as an optimisation problem with error function
$$ E({\asvector \Theta}) = \frac{1}{s} \sum_{i=1}^s \dist^2({\asvector \Theta}, {\asvector x^{(i)}}).$$
In Euclidean space, $\frechet {\collection x}$ is the unique minimum of $E$.
For Hadamard manifolds, such as hyperbolic space, the function $E$ is convex \cite{Bacak}, and consequently, also in this case, $E$ has a unique minimum $\frechet {\collection x} = \argmin_\Theta E({\asvector \Theta})$.
In this context, $\frechet {\collection x}$ is known as the {\it Fréchet mean} of the points.
\end{subsection}

\begin{subsection}{Optimisation}
The parameter vector ${\asvector \Theta} \in \hyperboloid n$ is considered as a point in the hyperboloid model.
It follows immediately from \eqref{hyperbolicdistancegradient} that, for any ${\asvector \Theta} \in \hyperboloid n$, the gradient of $E({\asvector \Theta})$ in Minkowski space is 
$$ \gradient {\asvector \Theta} {\minkowski n} E = \frac{2}{s} \sum_{i=1}^s -\dist({\asvector \Theta}, {\asvector x^{(i)}}) \cdot \left({{\minkowskiformpower n {\asvector \Theta} {{\asvector x^{(i)}}} 2} - 1}\right)^{-1/2} \cdot {\asvector x^{(i)}}  $$
The hyperboloid gradient $\gradient {\asvector \Theta} {\hyperboloid n} E$ is obtained by the projection \eqref{hyperboloidtangentspaceprojection}, as before.
Given a learning rate $\alpha > 0$, a new parameter vector $\new {\asvector \Theta}$ can be obtained either via the exponential map
$$ \new {\asvector \Theta} = \Exp_{\asvector p} (-\alpha \cdot \gradient {\asvector \Theta} {\hyperboloid n} E), $$
or by employing the retraction updates.
Beginning at the same initial parameter vector ${\asvector \Theta}^{(0)}$, these two distinct methods will yield separate sequences of parameter vectors.

The retraction update is effected by passing to the Poincaré ball model, performing a retraction update using \eqref{retraction}, and then returning again to the hyperboloid.
That is, the current parameter vector ${\asvector \Theta} \in \hyperboloid n$ and the gradient $\gradient {\asvector \Theta} {\hyperboloid n} E \in \tangentspace p \hyperboloid n$ are mapped to the Poincaré disc via \eqref{hyperboloidtopb} and \eqref{differentialhyperboloidtopb}, yielding a parameter vector ${\asvector \Psi} \in \pb n$ and a tangent vector at that point:
$$ {\asvector \Psi} = \projection ( {\asvector \Theta} ), \qquad \gradient {\asvector \Psi} {\pb n} {\tilde E} = \differential \projection |_{\asvector \Theta} \left(\gradient {\asvector \Theta} {\hyperboloid n} E\right),$$
where $\tilde{E} = E \circ \projection^{-1}$ is the error function considered on the Poincaré ball.
The retraction update is then calculated via
$$ \new {\asvector \Psi} = {\asvector \Psi} - \alpha \gradient {\asvector \Psi} {\pb n} {\tilde E}$$
where, if $\new {\asvector \Psi} \notin \pb n$ is no longer in the Poincaré ball, it is rescaled to have length $1 - 10^{-5}$, as per \cite{NickelKiela}.
The new parameter vector on the hyperboloid is then $\new {\asvector \Theta} = \projection^{-1} (\new {\asvector \Psi})$, calculated via \eqref{pbtohyperboloid}.
\end{subsection}

\begin{subsection}{Uniform sampling in hyperbolic space}
As observed earlier, the expected error of a retraction update depends upon the location of the point being updated.
In order to compare the two update methods, it is therefore necessary to consider the optimisation in different regions of hyperbolic space.
This is achieved by sampling a point $\asvector p$ on the hyperboloid from the uniform distribution on a ball centred at the hyperboloid base point with radius $\maxsubscript r$, and then drawing the samples ${\asvector x^{(i)}} \in \collection x$ from the uniform distribution on the ball with the same radius centred at $\asvector p$.
We restrict ourselves to the case where $n=2$, for simplicity.
The centre point $\asvector p \in \hyperboloid 2$ is sampled 50 times, and for each $\asvector p$, 50 collections $\collection x$ are constructed, each time by sampling $s$ times from the uniform distribution on the disc centred at $\asvector p$.
In view of the boundness of the error of the retraction update (c.f.\ Figure \ref{retractionerrorfigure}), we take $\maxsubscript r = 3$.

Uniform sampling on a disc may be achieved by first sampling a direction emanating from the centre point, and then sampling a distance $R \in [0, \maxsubscript r ]$ from the distribution on that interval whose cumulative distribution function (CDF) gives the proportion of the maximal disc area covered by a concentric disc of a given radius $r$.
That is, such that
$$p := \mathrm{P}(R < r) = \frac{\area r}{\area {\maxsubscript r}}, \quad r \in [0, \maxsubscript r ],$$
where $\area r$ is the area of a disc of radius $r$.
In the case of hyperbolic plane ($n=2$), we have $\area r = 2 \pi (\cosh r - 1)$, and so
\begin{equation}\label{cdf}
	p = \frac{\cosh r - 1}{\cosh {\maxsubscript r} - 1}.
\end{equation}
The distance $r$ is sampled via {\it inversion sampling}, that is, by solving \eqref{cdf} for $r$, yielding
\begin{equation}\label{cdftransposed}
	r = \arccosh(1 + p (\cosh(\maxsubscript r) - 1)).
\end{equation}
Samples of $r$ from the distribution with the CDF \eqref{cdf} can then be obtained from uniform samples of $p$ from the interval $[0,1]$, transformed to values of $r$ via \eqref{cdftransposed}.
\end{subsection}

\begin{subsection}{Steps until arrival}
We measure the expected number of steps to arrive within a neighbourhood of the solution.
Firstly, we compute the solution using the exponential updates, using a low learning rate and as many steps as necessary to converge (where convergence is detected by the vanishing of the gradient).
For each choice of centre point $\asvector p \in \hyperboloid n$, and each collection of samples $\collection x$, and for each of the two update methods (exponential or retraction), we then count the number of steps required for the optimisation to arrive within a distance of $10^{-4}$ of the solution, using a constant learning rate $\alpha$.
Learning rates $\alpha$ were swept over a range of values beginning at $\alpha=0.2$ and ending when $\alpha$ was sufficiently high that neither method could reliably arrive in the neighbourhood of the solution.
\end{subsection}

\begin{subsection}{Results}
Table \ref{table:mean_ttas} shows the mean number of steps for the two update methods at different values of the learning rate $\alpha$.
The optimal learning rates for the exponential- and retraction- updates were $\alpha=0.7$ and $\alpha=0.6$, respectively, and at their respective optimal learning rates, the estimated number of steps was $7.2$ and $12.8$.
Thus, on average, the exponential updates arrive more quickly than the retraction updates.
Figure \ref{stepstilarrivalscatterplots} provides a closer look at the comparative performance of the two methods at these two learning rates.
Even considering the learning rate $\alpha=0.6$ that is chosen for optimality of the retraction updates, we see that in 95.5\% of the trials the exponential updates arrive in the neighbourhood of the solution first.
From the slope of the line of best fit, it is apparent that approximately 46\% less updates are required if using the exponential map instead of the retraction.

\begin{table*}[t!]
\centering
\begin{tabular}{ c  c  c  c  c  c  c  c  c  c }
\toprule
Method / $\alpha$ & 0.20 & 0.30 & 0.40 & 0.50 & 0.60 & 0.70 & 0.80 & 0.90 & 1.00\\
\midrule
exponential & 34.4 & 21.8 & 15.1 & 10.2 & 7.4 & \textbf{7.2} & 7.9 & inf & inf \\
retraction & 35.4 & 22.8 & 16.7 & 13.7 & \textbf{12.8} & 15.3 & inf & inf & inf \\
\bottomrule
\end{tabular}
\caption{Mean number of updates before arriving within distance $0.0001$ of the Fréchet mean for the exponential and retraction updates.  The best performance for each method is marked bold.  The $5$ points were sampled uniformly from a disc of radius $3$.}
\label{table:mean_ttas}
\end{table*}

\begin{figure}
	\centering
	\includegraphics[width=7cm]{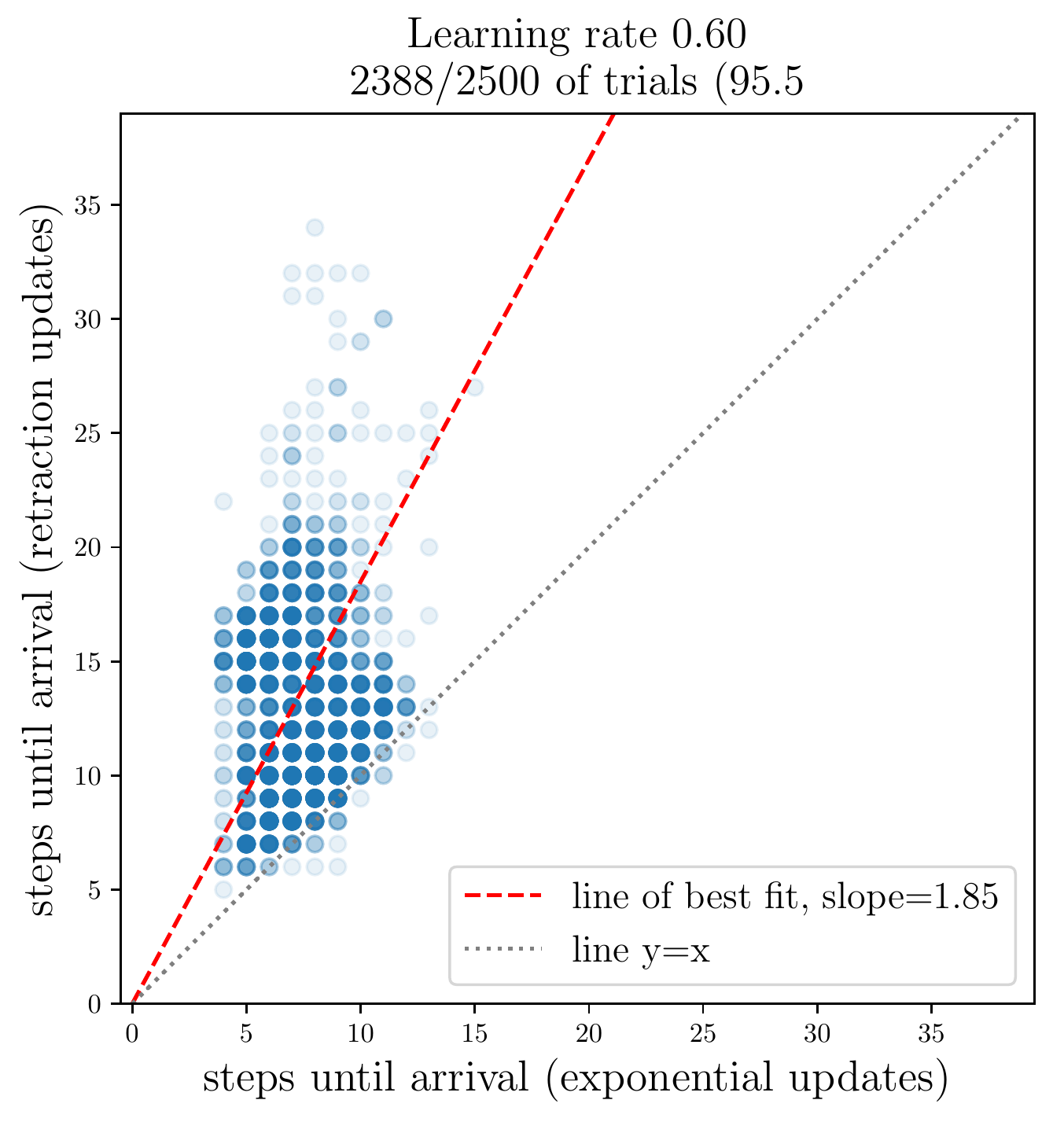}
    	\includegraphics[width=7cm]{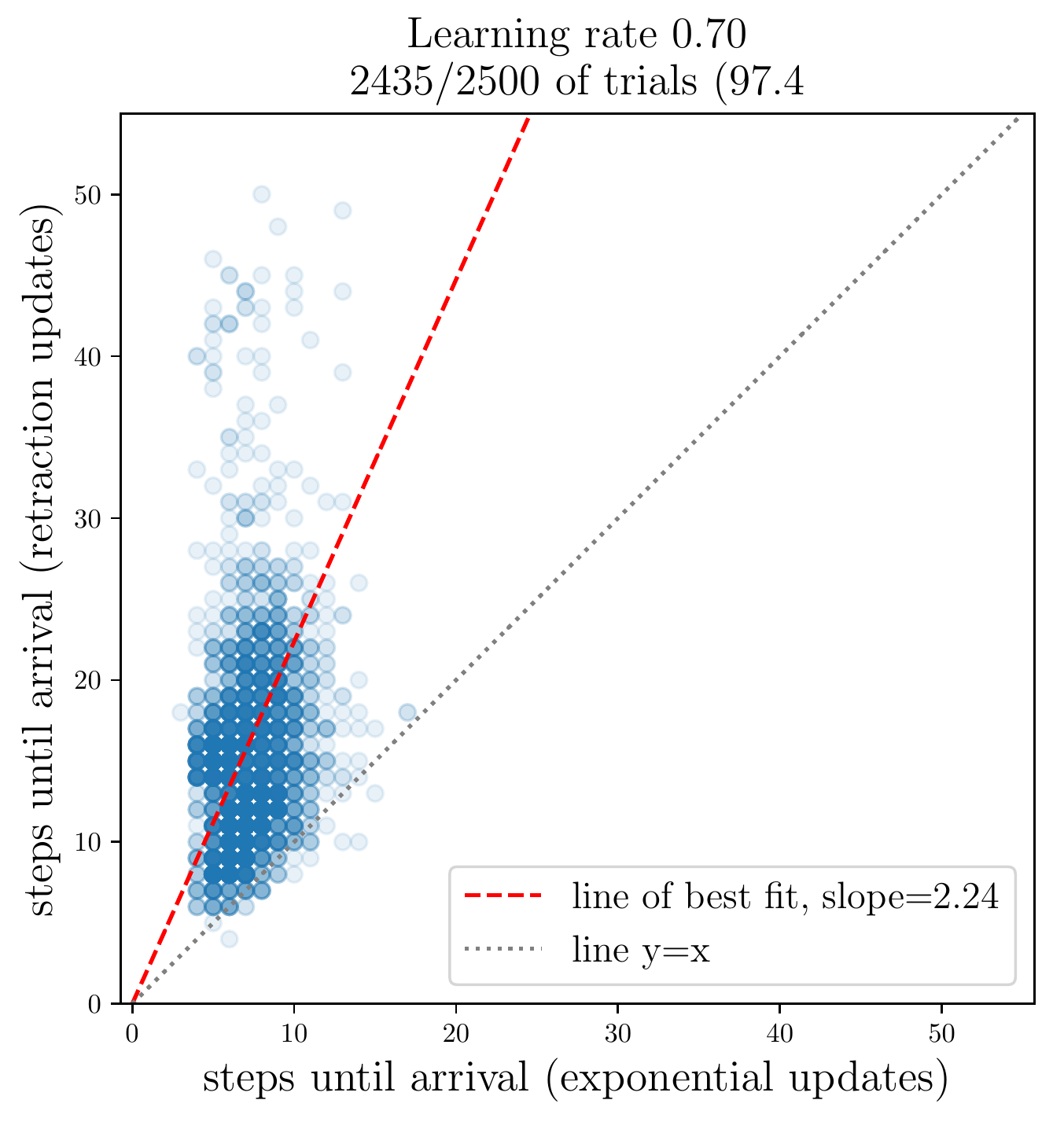}
	\captionof{Figure}{Scatter plots of the number of updates before arriving within distance $0.0001$ of the Fréchet mean for the exponential and retraction updates.  The $5$ points were sampled uniformly from a disc of radius $3$. Learning rates chosen for the optimality of the two methods from Figure \ref{table:mean_ttas}.
		Darkness of a point indicates its multiplicity.
		In the vast majority of cases, the exponential updates arrival in the neighbourhood of the solution first.
		For example, for learning rate $0.60$, the exponential updates arrived before the retraction updates in 95.5\% of the trials and the exponential updates arrived in approximately 46\% less steps (calculated from the slope). \label{stepstilarrivalscatterplots}}
\end{figure}

\end{subsection}
\end{section}

\begin{section}{Discussion and outlook}
We've seen that the hyperboloid model of hyperbolic space is very convenient for optimisation.
Specifically, it has been demonstrated that the calculations of distance and of the distance gradient are straight-forward in the hyperboloid model, and moreover that gradient descent (i.e.\ updating via the exponential map) is easy to implement there.
In particular, the hyperboloid model is much more appropriate for optimisation than the Poincaré ball model, though this has been preferred in some past optimisations.

It would certainly be of interest to implement the exponential updates in existing optimisations in hyperbolic space and to measure the improvement in performance.
The uniform sampling method described in section \ref{section:numericalexperiments} could also be useful for initialisation of model parameters, though further work is required to generalise the method to higher dimensions (perhaps using rejection sampling). 
\end{section}

\bibliographystyle{alpha}
\bibliography{references}

\end{document}